\input amstex
\input Amstex-document.sty

\pageno 293

\def\shorttitle{Differential Geometry via Harmonic Functions}
\def\shortauthor{P. Li}

\topmatter
\title\nofrills{\boldHuge Differential Geometry via Harmonic Functions}
\endtitle

\author \Large P. Li* \endauthor

\thanks *Department of Mathematics, University of
California, Irvine, CA 92697-3875, USA. E-mail: pli\@math.uci.edu \endthanks

\abstract\nofrills \centerline{\boldnormal Abstract}

\vskip 4.5mm

In this talk, I will discuss the use of harmonic functions to study the geometry and topology of complete manifolds.  In my previous joint work with Luen-fai Tam, we discovered that the number of infinities of a complete manifold can be estimated by the dimension of a certain space of harmonic functions.  Applying this to a complete manifold whose Ricci curvature is almost non-negative, we showed that the manifold must have finitely many ends.  In my recent joint works with Jiaping Wang, we successfully applied this general method to two other classes of complete manifolds.  The first class are manifolds with the lower bound of the spectrum $\lambda_1(M) >0$ and whose Ricci curvature is bounded by
$$
Ric_M \ge -{m-2 \over m-1} \lambda_1(M).
$$
The second class are stable minimal hypersurfaces in a complete manifold with non-negative sectional curvature.  In both cases we proved some splitting type theorems and also some finiteness theorems.

\vskip 4.5mm


\endabstract
\endtopmatter

\document

\baselineskip 4.5mm \parindent 8mm

\specialhead \noindent \boldLARGE 1. Introduction \endspecialhead

In 1992, the author and Luen-fai Tam \cite{12} discovered a general method to determine if a complete, non-compact, Riemannian manifold have finitely many ends.  An end is simply defined to be an unbounded component of the compliment of a compact set in the manifold.  If the number of ends is finite, their technique also provides an estimate on the number of ends.  In particular, they applied this method to prove that a certain class of manifolds must have finitely many ends.

\proclaim{Theorem 1 (Li-Tam)} Let $M^m$ be a complete, non-compact, manifold with
$$
Ric_M(x) \ge -k(r(x)),
$$
where $k(r)$ is a continuous non-increasing function satisfying
$$
\int_0^\infty r^{m-1}\,k(r)\, dr < \infty.
$$
Then there exists a constant  $0< C(m, k) < \infty$  depending only on $m$ and $k$, such that,
$M$ has at most $C(m, k)$ number of ends.
\endproclaim

Since a manifold with non-negative Ricci curvature will satisfy the hypothesis, this theorem can be viewed as a perturbed version of the splitting theorem \cite{4} of Cheeger-Gromoll.  A weaker version of the above theorem for manifolds with non-negative Ricci curvature outside a compact set was also independently proved by Cai \cite{1}.

In some recent work of Jiaping Wang and the author, they successfully applied the general theory of determining the number of ends to other situations.  The purpose of this note is to give a quick overview of the theory and its applications to manifolds with positive spectrum and minimal hypersurfaces.

\specialhead \noindent \boldLARGE 2. General theory
\endspecialhead

Throughout this article, we will assume that
$(M^m,\  ds_M^2)$ is an $m$-dimensional,  complete, non-compact
Riemannian manifold without boundary.
In terms of local coordinates $(x_1, x_2, \dots, x_m)$, if the metric is given by
$$
ds_M^2 = g_{ij}\,dx_i\,dx_j,
$$
then the  Laplacian  is defined by
$$
\Delta = {1\over\sqrt{g}} \,{\partial\over {\partial x_i}} \bigg( g_{ij}\,\sqrt{g} {\partial
\over {\partial x_j}}\bigg),
$$
where $(g^{ij}) = (g_{ij})^{-1}$ and $g = \det (g_{ij}).$
A function is said to be {\it harmonic}   on $M$ if it satisfies the Laplace equation
$$
\Delta f(x) = 0
$$
for all $x \in M.$

In order to state the general theorem, it is necessary for us to define the following spaces.
\proclaim{Definition 1} Let
$$
\Cal H_D(M)  = \{ f \, |\, \Delta f = 0, \, ||f||_\infty < \infty, \int_M |\nabla f|^2 < \infty\}
$$
be the space of bounded harmonic functions with finite Dirichlet integral defined on $M.$
\endproclaim

\proclaim{Definition 2} Let
$$
\Cal H_+(M)  =  \langle\{ f \, |\, \Delta f = 0, \, f >0\}\rangle
$$
be the space spanned by the set of positive harmonic functions defined on $M$.
\endproclaim

\proclaim{Definition 3} Let
$$
\Cal H'(M) = \langle\{ f \, |\, \Delta f = 0, \text{bounded from
one side on each end} \}\rangle
$$
be the space spanned by the set of harmonic functions defined on $M$, which has the property that each one is bounded either from above or below on each end.
\endproclaim

Observe that these spaces are monotonically contained in each other, i.e.,
$$
\Cal H_D(M) \subset \Cal H_+(M) \subset \Cal H'(M).
$$
Let us also recalled the following potential theoretic definition.
\proclaim{Definition 4} An end $E$ of $M$ is  {\it non-parabolic}  if it admits a positive Green's function with
Neumann boundary condition on $\partial E$.
Otherwise, $E$ is said to be  {\it parabolic}.
\endproclaim

It is important to note that if $M$ has at least one non-parabolic end,
 then $M$ admits a positive Green's function.  In this case, we say that $M$ is non-parabolic.  The interested reader can refer to \cite{11} for more detail descriptions.  Let us now state the general theorem in \cite{12}.

\proclaim{Theorem (Li-Tam)} Let $M$ be a complete, non-compact manifold without boundary.  Then there exists a subspace $\Cal K \subset \Cal H'(M)$, such that, $\dim \Cal K$ is equal to the number of ends of $M.$

Moreover, if $M$ is  non-parabolic , then the subspace $\Cal K$ can be taken to be in
$\Cal H_+(M)$.  Also there exists another subspace
$\Cal K_N \subset \Cal H_D(M)$, such that, $\dim \Cal K_N$
 is equal to the number of non-parabolic ends of $M$.
\endproclaim

At this point, it is important to point out that even though an estimate on the dimension of the spaces $\Cal H'(M)$, $\Cal H_+(M),$ or $\Cal H_D(M)$ will imply an estimate on the number of ends of corresponding type, however, in general, these spaces can be bigger than $\Cal K$ or $\Cal K_N.$  Hence to effectively use the above theorem, one should use the constructive argument in the proof of the theorem to give an estimate on $\Cal K$ and $\Cal K_N$ directly.  Indeed, this was the case in the proof of Theorem 1.  This is also true for the two applications stated in  the subsequence sections.

\specialhead \noindent \boldLARGE 3. Manifolds with positive
spectrum \endspecialhead

A complete manifold $(M, ds_M^2)$ is {\it conformally compact}  if $M$ is topologically a
manifold with boundary given by $\partial M.$ Moreover, there is a background metric
$ds^2_0$ on $(M, \partial M)$ such that
$$
ds_M^2 = \rho^{-2}\,ds^2_0,
$$
where $\rho$ is a defining function for $\partial M$ satisfying the conditions
$$
\rho = 0 \quad \text{on} \quad \partial M
$$
and
$$
 d \rho \neq 0 \quad  \text{on} \quad \partial M.
$$

A direct computation reveals that the sectional curvature, $K_M$, of the complete metric $ds^2$ has asymptotic value given by
$$
K_M \sim -|d \rho |^2,
$$
near $\partial M.$ Hence if $(M, ds_M^2)$ is also assumed to be Einstein with
$$
Ric_M = -(m-1),
$$
then
$$
K_M(x) \sim -1,
$$
as $x\to \infty.$

In 1999, Witten-Yau \cite{19} proved a theorem concerning the AdS/CFT correspondence, which effectively ruled out the existence of worm holes.  It is also a very interesting theorem in Riemannian geometry.

\proclaim{Theorem (Witten-Yau)} Let $M^m$ be a conformally compact, Einstein manifold of dimension at least 3.  Suppose the boundary $\partial M$ of $M$ has positive Yamabe constant,  then
$$
H_{m-1}(M, \Bbb Z) = 0.
$$
In particular, this implies that $\partial M$ is connected and $M$ must have only 1 end.
\endproclaim

Shortly after, Cai-Galloway \cite{2} relaxed the assumption of Witten-Yau by assuming the boundary $\partial M$ has non-negative Yamabe constant.  We would also like to point out that by a theorem of Schoen \cite{17}, a compact manifold has non-negative Yamabe constant is equivalent to the fact that it is conformally equivalent to a manifold with non-negative scalar curvature.

In his Stanford thesis, X. Wang \cite{18}  generalized the Witten-Yau, Cai-Galloway theorem by studying $L^2$ harmonic 1-forms.

\proclaim{Theorem (Wang)}
Let $M^m$ be a conformally compact manifold of dimension at least 3. Suppose the Ricci curvature of $M$  is bounded by
$$
Ric_M \ge -(m-1)
$$
and the lower bound of the spectrum of the Laplacian $\lambda_1(M)$ has a positive lower bound given by
$$
\lambda_1(M) \ge (m-2),
$$
then either

{\rm (1)} $M$  has no non-constant $L^2$-harmonic 1-forms, i.e.,
$$
H^1(L^2(M)) = 0;
$$
or

{\rm (2)} $M = \Bbb R \times N$    with the warped product metric
$$
ds_M^2 = dt^2 + \cosh^2 t\,ds_N^2,
$$
where $(N, ds_N^2)$  is  a compact  manifold with $Ric_N \ge -(m-2).$  Moreover, $\lambda_1(M) = m-2.$
\endproclaim

To see that this is indeed a generalization of the theorems of Witten-Yau and Cai-Galloway, one uses a theorem of Mazzeo \cite{16} asserting that on a conformally compact manifold
$$
H^1(L^2(M)) \simeq  H^1(M, \partial M).
$$
By a standard exact sequence argument, the conclusion  that $H^1(L^2(M)) = 0$ implies that $M$ has only 1 end.
In addition to this, one also uses a theorem of Lee \cite{10} giving a lower bound on $\lambda_1$ for conformally compact, Einstein manifold with non-negative Yamabe constant on $\partial M.$

\proclaim{Theorem (Lee)} Let $M$ be a conformally compact, Einstein manifold   with
$$
Ric_M = -(m-1).
$$
Suppose that $\partial M$ has  non-negative Yamabe constant, then
$$
\lambda_1(M) \ge {(m-1)^2 \over 4}.
$$
\endproclaim

Since ${(m-1)^2 \over 4} \ge m-2$, Wang's theorem implies the theorems of Witten-Yau and Cai-Galloway.  Observe that the warped product case in Wang's theorem has negative Yamabe constant on $\partial M$.

At this point, let us also recall a theorem of  Cheng \cite{5} stating that:
\proclaim{Theorem (Cheng)} Let $M$ be a complete manifold with
$$
Ric_M \ge -(m-1),
$$
then
$$
\lambda_1(M) \le {(m-1)^2 \over 4}.
$$
\endproclaim

Combining the results of Cheng and Lee
we conclude that
$$
\lambda_1(M) = {(m-1)^2 \over 4}
$$
for conformally compact, Einstein manifolds, whose Ricci curvature is given by
$$
Ric_M = -(m-1)
$$
and has non-negative Yamabe constant for its boundary.

In the authors recent joint work with Jiaping Wang \cite{14} , they  proved this splitting type theorem without assuming the manifold is conformally compact.

\proclaim{Theorem 2 (Li-Wang)} Let $M^m$ be a complete manifold with dimension $m\ge 3$.
Suppose the Ricci curvature of $M$ is bounded by
$$
Ric_M \ge -(m-1)
$$
and
$$
\lambda_1(M) \ge m-2,
$$
then either

{\rm (1)} $M$ has only 1 end with  infinite volume;

\noindent or

{\rm (2)}$M = \Bbb R \times N$    with the warped product metric
$$
ds_M^2 = dt^2 + \cosh^2 t\,ds_N^2,
$$
where $(N, ds_N^2)$  is  compact  with $Ric_N \ge -(m-2).$  Moreover, $\lambda_1(M) = m-2.$
\endproclaim

It is worth noting that this theorem implies that when the lower
bound for $\lambda_1(M)$ of Cheng is achieved, then either

{\rm (1)}$M$ has only  1 end with infinite volume,

\noindent or

{\rm (2)}$M= \Bbb R \times N$  is the  warped product   and
$m=3.$

Also, since all the ends of a conformally compact manifold must
have infinite volume, Theorem 2 is, in fact, a generalization of
the theorems of Witten-Yau, Cai-Galloway, and Wang.  It is also
interesting to note that without the conformally compactness
assumption, it is possible to have finite volume ends as indicated
by following example.
\proclaim{Example 1}\rm  Let $M^m = \Bbb R
\times N^{m-1}$ with the warped product metric
$$
ds_M^2 = dt^2 + \exp(2t)\,ds_N^2,
$$
where $N$ is a compact manifold with
$$
Ric_N \ge 0.
$$
\endproclaim

A direct computation shows that $M$ has Ricci curvature bounded by
$$
Ric_M \ge -(m-1)
$$
and
$$
\lambda_1(M) \ge m-2.
$$
In fact, when $m=3,$ $\lambda_1(M) = 1.$  Obviously $M$ has two ends. One end $E$ has infinite volume growth with
$$
V_E(r) \sim C\,\exp((m-1)\,r),
$$
while the other end $e$ has finite volume with volume decay given by
$$
V_e(\infty) - V_e(r) \sim C\, \exp(-(m-1)\,r).
$$

We would like to point out that the pair of conditions
$$
\left. \aligned Ric_M &\ge -(m-1)\\
\lambda_1(M) &\ge m-2
\endaligned \qquad \right\}\tag 1
$$
is equivalent to the pair of conditions
$$
\left. \aligned Ric_M &\ge -{m-1 \over m-2}\,\lambda_1(M)\\
\lambda_1(M) &>0.
\endaligned \qquad \right\} \tag 2
$$

On the other hand, the pair of conditions
$$
\left. \aligned Ric_M &\ge -{m-1 \over m-2}\,\lambda_1(M)\\
\lambda_1(M)&=0
\endaligned \qquad \right\} \tag 3
$$
are equivalent to the single assumption that
$$
Ric_M \ge 0,
$$
because the condition $\lambda_1(M) = 0$
is a consequence of the curvature assumption.

Taking this point of view, Theorem 2 can be viewed as an analogue to the splitting theorem of Cheeger-Gromoll.  Similarly to the fact that Theorem 1 is a perturbed version of the Cheeger-Gromoll splitting theorem, the following theorem in \cite{14} is a perturbed version of Theorem 2.

\proclaim{Theorem 3 (Li-Wang)} Let $M^m$ be a complete manifold with $m \ge 3.$ Suppose $B_p(R) \subset M$ is a geodesic ball such that
$$
Ric_M \ge -(m-1) \qquad    \text{on} \qquad  M\setminus B_p(R)
$$
and the lower bound of the spectrum of the Dirichlet Laplacian on $M\setminus B_p(R)$ is bounded by
$$
\lambda_1(M\setminus B_p(R)) \ge m-2 + \epsilon
$$
for some $\epsilon>0.$
Then there exists a constant $ 0< C(m, R, \alpha, v, \epsilon) < \infty $ depending only on $m$, $R$, $\alpha = \inf_{B_p(3R)} Ric_M,$ $ v = \inf_{x\in B_p(2R)} V_x(R)$, and $\epsilon$, so that the number of infinite volume ends of $M$ is at most $C(m, R, \alpha, v, \epsilon)$.
\endproclaim

In both Theorem 2 and Theorem 3, the authors only managed to estimate the number of infinite volume ends by estimating the number of non-parabolic ends.  In fact, when a manifold has positive spectrum, they proved that an end must either be non-parabolic with exponential volume growth, or it must be parabolic and finite volume with exponential volume decay.  Moreover, these growth and decay estimates can be localized at each end.

\proclaim{Theorem 4 (Li-Wang)} Let $M$ be a complete, non-compact,
Riemannian manifold.  Suppose $E$ is an end of $M$ given by a
unbounded component of $M \setminus B_p(R)$, where $B_p(R)$ is a
geodesic ball of radius $R$ centered at some fixed point $p \in
M$.  Assume that the lower bound of the spectrum $\lambda_1(E)$ of
the Dirichlet Laplacian on $E$ is positive.  Then as $r \to
\infty,$ either

{\rm (1)} $E$ is non-parabolic and has volume growth given by
$$
V_E(r) \ge C_1\, \exp (2\sqrt{\lambda_1(E)}\, r)
$$
for some constant $C_1>0$;

\noindent or

{\rm (2)} $E$ is parabolic and has finite volume with volume decay
given by
$$
V(E) - V_E(r) \le C_2\, \exp(-2\sqrt{\lambda_1(E)}\, r)
$$
for some constant $C_2 >0.$

In particular, if $\lambda_1(M) >0$, then $M$ must have exponential volume growth given by
$$
V_p(r) \ge C_1\, \exp(2\sqrt{\lambda_1(M)}\, r).
$$
\endproclaim

Both the volume growth and the volume decay estimates are sharp.  For example, the growth estimate is achieved by the hyperbolic $m$-space, $\Bbb H^m.$  Also, in Example 1 when dimension $m=3$, the infinite volume end achieves the sharp volume growth estimate and the finite volume end achieves the sharp volume decay estimate.  It is also interesting to point out that the sharp volume growth estimate is previously not known for manifolds with $\lambda_1(M)>0.$

 \specialhead \noindent \boldLARGE 4. Minimal hypersurfaces \endspecialhead
Let us recall that the well-known Bernstein's theorem (Bernstein, Fleming, Almgren, DeGiorgi, Simons) asserts that an entire minimal graph  $M^m \subset \Bbb R^{m+1}$ must be linear if $m\le 7.$  Moreover, the dimension restriction is necessary as indicated by the examples of Bombieri, DeGiorgi, and Guisti.
Since minimal graphs are necessarily area minimizing and hence stable (second variation of the area functional is non-negative), Fischer-Colbrie and Schoen \cite{8} considered a generalization of Bernstein's theorem in this category.  They proved that a  complete, oriented, immersed, stable minimal surface  in  a complete manifold with  non-negative scalar curvature must be conformally equivalent to either $\Bbb C$ or $\Bbb R \times \Bbb S^1.$
Moreover, if the ambient manifold is $\Bbb R^3$ then the minimal surface must be planar. This special case  was independently proved by
do Carmo and Peng \cite{6}.

Later, Fischer-Colbrie \cite{7} studied the structure of minimal surfaces with finite index.  Recall that a minimal surface has finite index means that there are only a finite dimension of variations such that the second variations of the area functional is negative.  In this case, Fischer-Colbrie proved that a  complete, oriented, immersed, minimal surface with finite index in a complete manifold with  non-negative scalar curvature must be conformally equivalent to a compact Riemann surface with finitely many punctures.  In particular, $M$ must have finitely many ends. The special case when $N= \Bbb R^3$ was also independently proved by
Gulliver \cite{9}.     It is in the spirit of the number of ends that Cao, Shen and Zhu \cite{3} found a  higher dimensional statement for stable minimal hypersurfaces in $\Bbb R^{m+1}.$

\proclaim{Theorem (Cao-Shen-Zhu)} Let
$M^m \subset \Bbb R^{m+1}$ be a complete, oriented,  immersed,  stable minimal
hypersurface in $\Bbb R^{m+1},$ then $M$ must have only 1 end.
\endproclaim

This theorem is recently generalized to minimal hypersurfaces with finite index by the author and Jiaping Wang \cite{13}.

\proclaim{Theorem 5 (Li-Wang)} Let
$M^m \subset \Bbb R^{m+1}$ be a complete, oriented,  immersed,  minimal
hypersurface with finite index in $\Bbb R^{m+1}$, then  $M$ must have finitely many ends.
\endproclaim

In another paper \cite{15}, they also considered complete, properly immersed, stable (or with finite index) minimal hypersurfaces in a complete, non-negatively curved manifold.

\proclaim{Theorem 6 (Li-Wang)} Let
$M^m \subset  N^{m+1}$ be a complete, oriented,
 properly immersed,  stable, minimal hypersurface.  Suppose
$N$ is a  complete manifold with  non-negative sectional
curvature.  Then either

{\rm (1)}  $M$ has only 1 end;

\noindent  or

{\rm (2)} $M=\Bbb R \times S$   with the product metric, where $S$
is a compact manifold with non-negative sectional curvature.
Moreover, $M$ is  totally geodesic in $N.$
\endproclaim

\proclaim{Theorem 7 (Li-Wang)} Let
$M^m \subset  N^{m+1}$ be a complete, oriented,
 properly immersed,  minimal hypersurface with finite index.  Suppose
$N$ is a  complete manifold with  non-negative sectional
curvature.  Then $M$ must have finitely many ends.
\endproclaim

It is interesting to point out that in the case when $M = \Bbb R \times S,$ the manifold is parabolic.  In this case, it is necessary to estimate the space $\Cal K$ rather than $\Cal K'.$  Again, the crucial point is to follow the construction of $\Cal K$ and obtain sufficient estimates on the functions in $\Cal K$ so that analytic techniques can be applied.  In the case of Theorem 5, since the ambient manifold is $\Bbb R^{m+1}$ and hence the ends of $M$ must all be non-parabolic, it is sufficient to estimate the space $\Cal K'$ as stated in Theorem 2.

\specialhead \noindent \boldLARGE References \endspecialhead

\widestnumber\key{AA}

\ref \key{1} \by M. Cai
\paper \rm Ends of Riemannian manifolds with nonnegative Ricci curvature outside a compact set
\jour {\it Bull. AMS} \vol 24 \yr 1991 \pages 371--377
\endref
\ref\key{2} \by M. Cai and G. J. Galloway
\paper \rm Boundaries of zero scalar curvature in the ADS/CFT correspondence
\jour {\it Adv. Theor. Math. Phys.} \vol 3 \yr 1999 \pages 1769--1783
\endref

\ref\key{3} \by H. Cao, Y. Shen, and S. Zhu
\paper \rm The structure of stable minimal hypersurfaces in $\Bbb R^{n+1}$
\jour {\it Math. Res. Let.} \vol 4 \yr 1997 \pages 637--644
\endref

\ref \key{4} \by J. Cheeger and D. Gromoll
\paper \rm The splitting theorem for manifolds of nonnegative Ricci curvature
\jour {\it J. Diff. Geom.} \vol 6  \yr 1971 \pages 119--128
\endref

\ref \key{5} \by S. Y. Cheng
\paper  \rm Eigenvalue Comparison theorems and its Geometric Application
\jour {\it Math. Z.} \vol 143 \yr1975 \pages  289--297
\endref

\ref\key{6} \by M. do Carmo and C. K. Peng
\paper \rm Stable complete minimal surfaces in $\Bbb R^3$ are planes
\jour {\it Bull. AMS} \vol 1 \yr 1979 \pages 903--906
\endref

\ref \key{7} \by D. Fischer-Colbrie
\paper \rm On complete minimal surfaces with finite Morse index in three manifolds
\jour {\it Invent. Math.} \vol 82 \yr 1985 \pages 121--132
\endref

\ref \key{8} \by D. Fischer-Colbrie and R. Schoen \paper \rm The
structure of complete stable minimal surfaces in 3-manifolds of
non-negative scalar curvature \jour {\it Comm. Pure Appl. Math.}
\vol 33 \yr 1980 \pages 199--211
\endref

\ref\key{9} \by R. Gulliver
\paper \rm Index and total curvature of complete minimal surfaces.
\inbook {\it Geometric measure theory and the calculus of variations (Arcata, Calif., 1984), Proc. Sympos. Pure Math. 44}
\publ Amer. Math. Soc., \publaddr Providence, RI. \yr 1986 \pages 207--211
\endref
\ref\key{10} \by J. Lee
\paper \rm The spectrum of an asymptotic hyperbolic Einstein Manifold
\jour {\it Comm. Anal. Geom.} \vol 3 \yr 1995 \pages 253--271
\endref

\ref \key{11} \by P. Li \paper \rm Curvature and function theory on Riemannian manifolds \inbook {\it Survey in
Differential Geometry ``In Honor of Atiyah, Bott, Hirzebruch, and Singer"} \publ International Press \publaddr
Cambridge \vol VII \yr 2000, 71--111
\endref

\ref \key{12} \by P. Li and L. F. Tam
\paper \rm Harmonic functions and the structure of complete manifolds
\jour {\it J. Diff. Geom.} \vol 35 \yr 1992 \pages 359--383
\endref

\ref\key{13} \by P. Li and J. Wang
\paper \rm Minimal hypersurfaces with finite index
\jour {\it Math. Res. Let.} \vol 9 \yr 2002 \pages 95--103
\endref

\ref\key{14} \by P. Li and J. Wang
\paper \rm Complete manifolds with positive spectrum
\jour {\it J. Diff. Geom.} \vol 58 \yr 2001 \pages 501--534
\endref

\ref\key{15} \by P. Li and J. Wang
\paper \rm Stable minimal hypersurfaces in a nonnegatively curved manifold
\jour {\it Preprint}
\endref

\ref\key{16} \by R. Mazzeo
\paper \rm The Hodge cohomology of a conformally compact metric
\jour {\it J. Diff. Geom.} \vol 28 \yr 1988 \pages 309--339
\endref

\ref \key{17} \by R. Schoen
\paper \rm Conformal deformation of a Riemannian metric to constant scalar curvature
\jour {\it J. Diff. Geom.} \vol 20 \yr 1984 \pages 479--495
\endref

\ref\key{18} \by X. Wang
\paper \rm On conformally compact Einstein manifolds
\jour {\it Math. Res. Let.} \vol 8 \yr 2001 \pages 671--688
\endref

\ref\key{19} \by E. Witten and S. T. Yau
\paper \rm Connectedness of the boundary in the AdS/CFT correspondence
\jour {\it Adv. Theor. Math. Phys.} \vol 3 \yr 1999 \pages 1635--1655
\endref

\enddocument